\documentclass{article}
\usepackage{amsmath}
\usepackage{amsmath,amscd}
\usepackage{amssymb,amsfonts,amsmath,amscd}
\usepackage[margin=0.5in,footskip=0.25in]{geometry}
\usepackage{fancybox}
\usepackage{pgf,pgfarrows,pgfnodes,pgfautomata,pgfheaps,pgfshade}
\usepackage[all,graph,frame]{xy}
\usepackage{appendix}
\usepackage{enumerate}

\renewcommand{\theenumii}{\alph{enumii}}
\renewcommand{\theenumiii}{(\roman{enumiii})}

\def\R{\mathbb{R}}

\def\qed{{$\Box$}}

\newtheorem{theorem}{Theorem}

\newtheorem{conjecture}[theorem]{Conjecture}

\newtheorem{proposition}{Proposition}
\newtheorem{remark}{Remark}

\begin{document}


\title{Backward induction in presence of cycles}
\author{Vladimir Gurvich
\footnote
{National Research University, Higher School of Economics (HSE), Moscow Russia
             \textit{vgurvich@hse.ru}}}
\maketitle


\begin{abstract}
For the classical backward induction algorithm,
the input is an arbitrary $n$-person
positional game with perfect information
modeled by a finite acyclic directed graph
(digraph) and the output is a profile $(x_1, \ldots, x_n)$
of pure positional strategies that form
some special subgame perfect Nash equilibrium.
We extend this algorithm to work with
digraphs that may have directed cycles.
Each digraph admits a unique partition
into strongly connected components, which
will be treated as the outcomes of the game.
Such a game will be called
a {\em deterministic graphical multistage}(DGMS) game.
If we identify the outcomes corresponding
to all strongly connected components,
except terminal positions, we obtain
the so-called {\em deterministic graphical}(DG) games,
which are frequent in the literature.
The  outcomes of a DG game are all terminal positions
and one special outcome  $c$  that is assigned to all infinite plays.
We modify the backward induction procedure to adapt it for the DGMS games.
However, by doing so, we lose two important properties:
the modified algorithm always outputs a {\em Nash equilibrium}
(NE) only when  $n = 2$ and, even in this case,
this NE  may be not {\em subgame perfect}.
(Yet, in the zero-sum case it is.)
The lack of these two properties is not a fault of the algorithm,
just (subgame perfect) Nash equilibria
in pure positional strategies may fail to exist in the considered game.

\end{abstract}

\section{Introduction}
\label{s0}

The most natural think to do in the presence of potentially infinite play
seems to be to simply permit it...

Alan Washburn \cite{Was90}

and, even better, to make a typical play infinite.

\subsection{Decomposing digraphs into strongly connected components}
\label{ss00}
Let  $G = (V, E)$  be a finite {\em directed graph} (digraph) that
may have parallel edges, directed oppositely or similarly,
and loops,  but at most one at each vertex.
Digraph  $G$  is called {\em strongly connected}  if
for any  $v, v' \in V$  there is a directed path from  $v$  to  $v'$
(and, hence, from  $v'$ to $v$, as well).
By this definition, the union of two strongly connected digraphs
is strongly connected  whenever they have a common vertex.
A vertex-inclusion-maximal strongly connected induced subgraph of $G$
is called its {\em strongly connected component}.
Obviously, any digraph  $G = (V, E)$  admits a unique decomposition into such components:
$G^j = G[V^j] = (V^j, E^j)$  for $j \in J$, where
$J$ is a set of indices and
$V = \cup_{j \in J} V^j$  is a partition of  $V$, that is,
$V^j \cap V^{j'} = \emptyset$  whenever  $j \neq j'$  for a pair  $j, j' \in J$.
This partition can be determined in time linear in
the size of  $G$ (that is, in $(|V| + |E|)$) and
has numerous applications; see \cite{Tar72} and \cite{Sha81} for more details.
Here we suggest one more application, to positional games.

\medskip

Let  $V_T \subseteq V$  denote the set of all terminal vertices
(or terminals, for short); more precisely,
$v \in V_T$  if and only if there is no edge from  $v$,
except, perhaps, a loop  $e_v = (v, v)$.
Any such vertex forms a maximal connected component, that is,
$\{v\} = V^j$  for some  $j \in J$;
denote by  $J_T$  the corresponding subset of  $J$.
For convenience, let us add a loop  $e_v$  to each  $v \in V_T$
that does not have one yet.
All loops in  $V_T$,  old and new, will be called {\em terminal} loops
and all  $v \in V_T$  will be still called terminals.

The concepts of a {\em directed path}
and a {\em directed cycle} (dicycle) are defined standardly;
see, for example,  \cite{Tar72} and \cite{Sha81}.
An edge is a directed path of length $1$;
a loop is a dicycle of length  $1$,
a dicycle of length  $2$  is a pair of parallel oppositely directed edges,
a dicycle of length $|V|$  is called Hamiltonian,
A digraph without dicycles is called {\em acyclic}.

\smallskip

Let  $J_0$  denote the subset of  $J$  such that  $V^j$
is a single vertex with no loop it it.
Obviously, $G^j$ contains a dicycle unless  $j \in J_0$.
By our convention, $J_T \cap J_0 = \emptyset$.

For  each $j \in J$, let us contract  $G^j$
into a single vertex  $v^j$.
Then, all edges of  $E^j$
(in particular, the loops) disappear and
we obtain an acyclic digraph  $G^* = (V^*, E^*)$.

In case of the standard backward induction we have
$G^* = G$  or, equivalently,  $J = J_0 \cup J_T$  is
a partition of  $J$.

\subsection{Modeling positional games by digraphs}
\label{ss01}
Given an arbitrary digraph  $G = (V, E)$, we interpret
it as a {\em positional structure} in which
$V$  and  $E$  are positions and moves, respectively.
Furthermore, let us fix a partition
$V = V_1 \cup \ldots \cup V_n \cup V_T$  in which
$V_i$  are the sets of positions controlled by the players
$i \in I = \{1, \ldots, n\}$, while  $V_T$  is the set of all terminals.

The considered model is restricted to the games
with perfect information and without moves of chance.
Also, we restrict the players  $i \in I$
to their {\em pure and positional}
(also called {\em stationary} strategies  $X_i$.
Such a strategy  $x_i \in X_i$  assigns
a move $(v, v') \in E$  to every  $v \in V_i$.
Given an {\em initial position}  $v_0 \in V \setminus V_T$,
each {\em strategy profile}
$x = (x_1, \ldots, x_n)$  uniquely  defines
a walk in  $G$  that begins in  $v_0$  and
in each position follows the move chosen be the corresponding strategy.
Since all strategies are positional, this walk forms a ``lasso",
that is, it consists of an initial directed path and a dicycle repeated infinitely.
For every  $j \in J \setminus J_0$, let us identify all lassos whose dicycles
belong to  $G^j$ and treat them as a single outcome  $c_j$  of the game.
Then, $A = \{c_j \mid j \in  J \setminus  J_0\}$  is the set of outcomes.
An outcome is called a {\em  terminal} whenever $G^j$  is a terminal loop.
Thus, we obtain a {\em game form}  $g : \prod_{i \in I}  X_i  \rightarrow A$.
Positional structures, game forms, and games defined above
will be called {\em Deterministic Graphical Multistage} (DGMS) ones.

Let us identify the outcomes  $c_j$  for all  $j \in J \setminus J_T$
and denote the obtained new outcome  $c$;
it is assigned to each dicycle, or in other words,
to every infinite play, in  $G$.
As before, the terminal loops form distinct outcomes.
The corresponding positional structures, game forms, and games
are called {\em Deterministic Graphical} ones.
They were introduced by  Washburn \cite{Was90} in 1990
and since then are frequent in the literature; see, for example,
\cite{Was90}, \cite{Con92}, \cite{BG03}, \cite[Section 12]{BGMS07},
\cite{BG09}, \cite{AGH10}, \cite{AHMS12}, \cite{BEGM12}, \cite{Gur15}.
Only two-person zero-sum DG games were considered in
\cite{Was90}, \cite{Con92}, \cite{AHMS12}, while our
main results are related to the two-person,
but not necessarily zero-sum, DGMS games.

Note also that DG games are simpler than DGMS games,
since identification of any outcomes
respects the properties of solvability, which we study here.
Thus, ``DGMS games" could be an alternative title for this paper.

\subsection{Nash equillibria and solvability of game forms}
\label{ss02}
Given an arbitrary, not necessarily DG or DGMS, game form
$g : X \rightarrow A$, where $X = \prod_{i \in I}  X_i$,
let us introduce a utility function  $u : I \times A \rightarrow \R$.
Here the real number $u(i, o)$  is interpreted as a profit
of the player  $i \in I$  in case of the outcome  $a \in A$.
A pare  $(g, u)$  is called a {\em (normal form) game}.

A strategy profile  $x = (x_i \mid i \in I) \in X$
is called a {\em Nash equilibrium} (NE) if  $u(i, g(x)) \geq  u(i, g(x'))$
for every  $x' \in X$  that differs from  $x$  only in the coordinate $i$;
or in other words, if no player can strictly improve  by  replacing his/her strategy
provided all others keep their strategies unchanged.
This concept was suggested by Nash in  \cite{Nas50} and \cite{Nas51}.

If  $n = |I| = 2$  then  $g, u$,  and  $(g, u)$  are called
a two-person game form, utility function, and game, respectively.
Furthermore,  $u$  and  $(g, u)$   are called {\em zero-sum}
utility function and game, respectively, if
$u(1, a) + u(2, a) = 1$   for every outcome $a \in A$.
In this case a NE is usually referred to as a {{\em saddle point}.
Finally, $u$  and  $(g, u)$   are called {\em win/lose}
utility function and game, respectively,
if  $u$  is zero-sum and takes only values $+1$ and $-1$.

A game form  $g$  is called
\begin{itemize}
\item{$(S1)$} Nash-solvable (NS) if game  $(g,u)$  has a NE (in pure strategies) for any $u$;
\item{$(S2)$} zero-sum-solvable
if $n = 2$  and  $(g, u)$  has a saddle point for every zero-sum  $u$;
\item{$(S3)$} win/lose-solvable
if  $n = 2$ and  $(g, u)$  has a saddle point for every win/lose  $u$.
\end{itemize}


\begin{theorem}
\label{t1}
(\cite{Gur75}).
For the two-person game forms the above three properties are equivalent.
\end{theorem}

A two-person game form satisfying these properties
will be called just {\em solvable}, since no specification
NS, zero-sum, or win/lose is needed.

Implications $(S1) \Rightarrow (S2) \Rightarrow (S3)$
are immediate from the above definitions.
For the two-person zero-sum case, equivalence  $(S2)$ and $(S3)$
was  proven by Edmonds and Fulkerson
(who called it the ``Bottle-Neck Extrema Theorem") in \cite{EF70}  and later,
independently, in \cite{Gur73}.
Then implication $(S3) \Rightarrow (S1)$  was proven in \cite{Gur75}.
The same proof was repeated in a more focused paper \cite{Gur88}; see also \cite{RM-DS17}.

The proof is constructive:
a NE in the game  $(g, u)$ is obtained by an explicit algorithm
whenever the game form  $g$  is NS.
Unfortunately, this algorithm is exponential in the size of $A$.

\subsection{Main results. DG, DGMS, and other solvable game forms}
\label{ss03}
Another proof, based on a polynomial algorithm, appeared later;
see  \cite[Section 4.3]{Gur97}, \cite[Section 3]{BG03},
and the recent preprint \cite{GK17}.
Some new ideas from \cite{BBM89} and  \cite{DS91} were used in this proof.
For convenience of the reader, in Section \ref{ss11}
we will prove constructively a slightly stronger statement.

\begin{proposition}
\label{p1}
Implication  $(S3) \Rightarrow (S1)$  holds for the two-person game forms.
Moreover, to get a NE in a two-person game  $(g, u)$
with a solvable game form  $g$
and arbitrary  $u : \{1,2\} \times A \rightarrow \R$
we need to solve at most $2 |A|$  win/lose games
with the same game form  $g$.
\end{proposition}

Complexity of the latter problem depends on
the considered class of game forms.
We know several classes for which the win/lose games can be solved efficiently.

\medskip

Monotone bargaining was studied recently in \cite{GK17}, where
it was shown that all considered game forms are  NS  and
the corresponding  $\pm 1$ games are solvable in linear time
by a very simple greedy algorithm.

\medskip

Two-person positional game structures in which
distinct dicycles form distinct outcomes were considered in \cite{BGMS07}.
In this case only some special game forms are solvable.
They were characterized in \cite{BGMS07}, however,
not in general, but only in case of the so-called symmetric digraphs.
($G = (V, E)$  is called symmetric if
$(v, v') \in E$  whenever  $(v', v) \in E$).
Again the win/lose games are solvable in linear time \cite[Theorem 5]{BGMS07}.

\smallskip

We will derive solvability results for DGMS game forms
from the corresponding results for DG game forms.
The latter are known.
The two-person zero-sum DG games were introduced in \cite{Was90}.
For their pre-history, history, algorithms of solution, and their complexity,
we refer the reader to the short but comprehensive survey  \cite{AHMS12}.

Every win/lose DG game has a subgame perfect NE
(saddlle point) in pure  positional strategies, which
can be found in linear time in the size of  $G$,
that is, $|V| + |E|$.
This result was obtained in 1992 by Condon \cite{Con92}.
In fact, the same algorithm was discovered earlier, in 1986,
by Thompson \cite{Tho86} and presented to
the artificial intelligence community, under the name of
{\em retrograde  analysis}; see \cite{AHMS12} for more details..

Applying theorem \ref{t1}, we conclude that
the two-person DG game forms are NS, which is also known:
\cite[Section 3]{BG03}, \cite{BG09}, \cite[Section 12]{BGMS07}.
By backward induction, we will extend this result
from DG to DGMS case, and prove in Section \ref{ss12} the following statement.

\begin{proposition}
\label{p2}
The two-person DGMS game forms are NS and
the corresponding win/lose games can be solved in linear time.
\end{proposition}

Let us notice again that
the obtained NE may be not subgame perfect, but
in the zero-sum case it is.
To extend this result to the DGMS games, we use
some known results on the subgame perfect saddle points
of the two-person zero-sum DG games, which are of independent interest.

In fact, a DG game can be viewed as a
(very special case of a) mean-payoff stochastic game with perfect information.
Indeed, we can set the local rewards
$r(i, e) = 0$  for all  $i \in I$  and all  $e in E$,
except for the terminal loops  $e_v = (v, v)$,
for which the values of  $r(i, e_v)$ can be arbitrary.
Thus, the existence of a subgame perfect saddle point
in every zero-sum DG game follows from a much more general
(and old, 1957) result of Gillette \cite{Gil57}
claiming the same property for the stochastic games with perfect information.
(Gillette's proof, based on the Hardy-Littlewood Theorem \cite{HL31} of 1931,
had a flaw; conditions of the theorem were not accurately verified.
It was corrected in 1969 by Liggett and Lippman \cite{LL69}.)

\smallskip

An ``almost linear" algorithm was suggested in \cite[Sections 3] {AHMS12}.
It outputs a subgame perfect saddle point in time $|E| + |A|\log|A|)$ and requires
sorting of the outcomes, that is, $|A| \log |A|$  comparisons.
Given a permutation ordering the terminals, the problem becomes linear.

\smallskip

We extend the above results from DG to DGMS games,
combining the known DG algorithms with backward induction,
and prove in Section \ref{ss12} the following statement.

\begin{proposition}
\label{p3}
A subgame perfect saddle point in a two-person zero-sum
DGMS game can be found in time linear in  $(|E| + |V| + |A| \log |A|)$ and
in linear time if the outcomes are sorted in advance.
\end{proposition}

As we already mentioned,
one obtains a DGMS game from a DG one,
identifying the outcomes  $c_j$  for all
$j \in J \setminus (J_0 \cup J_T)$.
Obviously, any identification of outcomes respects solvability.

\subsection{Modifying Backward induction}
\label{ss04}
The backward induction algorithm was developed
for $n$-person games with perfect information
by Gale \cite{Gal53} and Kuhn \cite{Kuh53}; see also \cite{Kuh50}.
The authors model games by arborescence
(directed rooted trees) but an extension
to arbitrary acyclic digraphs is straightforward.
Moreover, positions of chance can be easily
included into consideration too \cite{Gal53}.
The algorithm constructs, in linear time,
some special subgame perfect NE in pure positional strategies.

However, the acyclicity is very essential.
In Section \ref{ss12} we modify the backward induction procedure
to allow it to work in presence of dicycles, but
for this we have to pay a price:
the modified backward induction works
only in the two-person case and the obtained NE
may be not subgame perfect.

First, we introduce the new procedure for the win/lose games.
Yet, we want to construct a NE for an arbitrary
(not necessarily zero-sum) two-person game,
and we will, but the algorithm becomes more sophisticated.
It calls at most $2 |A|$  times for solving some specially chosen win/lose games.
A constructive proof based on Theorem \ref{t1}  and Propositions \ref{p1} and \ref{p2}
will be given in Section \ref{s1}.
We always obtain a NE but it may be not subgame perfect,
although in the zero-sum case it is.

These above two flaws is not a fault of the algorithm.
Just NE may fail to exist when $n > 2$.
First such example, with  $n=4$, was constructed in \cite{Gur15};
then, a much smaller one, with  $n=3$, was given in \cite{BGMOV16};
see Section \ref{s2} for more details.
The main result of the present paper implies that
$n$  cannot be reduced further, that is,  NS holds for  $n = 2$.
Yet, a subgame perfect NE may fail to exist even in this case.
First such example was given in \cite[Section 6.2, Figure 10]{AGH10};
more complicated examples having some extra properties
can be found in  \cite{BEGM12} and \cite{GO14};
see Section \ref{ss05} for more details.

\subsection{Interpretation}
\label{ss05}
Members of a household  $I = \{1, \ldots, n\}$
have to delegate somebody to perform a necessary
(but, perhaps, not very pleasant) job, like
weekly shopping or cleaning the house.
We model the game by a digraph $G = (V, E)$  such that
$$V = \{v_1, \ldots, v_n; t_1, \ldots, t_n\}, \;
E = \{(v_1,v_2), \ldots, (v_{n-1},v_n), (v_n, v_1);
(v_1, t_1), \ldots, (v_n, t_n)\}.$$
We also add a terminal loop
$c_i$  in  $t_i$  for every  $i \in I$.
Digraph  $G$  has $n+1$  strongly connected components:
$n$  terminal loops and one $n$-dicycle.
Thus, $G^*$  is the directed star with  $n$  rays.
Each player $i \in I$  controls only one position, $v_i$, and
has only two possible moves in it;
in other words, (s)he has two strategies:
either to proceed along the $n$-dicycle or terminate in  $c_i$.
The latter means that (s)he accepts the job.
If no player does, then the job remains undone,
negotiations come to a dead-end and  the play,
in contrast, results in the long dicycle.
The corresponding outcome is  $c_0$  and
$A = \{c_0, c_1, \ldots, c_n\}$  is the set of all outcomes.
A utility function  $u : I  \times A \rightarrow \R$  may be arbitrary.
It would be natural  to assume that  $u(i, c_i) < u(i, c_j)$
for any distinct  $i, j \in \{1, \ldots, n\}$,
but even this is not required.

For  $n=2$  and  $n=3$  this example  was analyzed
in \cite[Section 6.2]{AGH10} and \cite{BG03}, respectively;
see also \cite{BEGM12} and \cite[Remark3]{BGMOV16}.
It was demonstrated that
a subgame perfect NE may fail to exist in both cases,
yet, with respect to any fixed initial position, it exists.
However, this property may fail too. 
First NE free example, with  $n=4$,  was constructed in \cite{Gur15};
then, a smaller one, with  $n=3$, in \cite{BGMOV16}.
The main result of the present paper implies that $n$ cannot be reduced further.

\medskip

In general, the acyclic digraph  $G^*$, as well as
the strongly connected components of  $G$  may have more complicated structure;
see, for example, \cite{BG03}.
Every strongly connected component  $G^j$
can be interpreted as a diplomatic issue.
(The level may not be very high, as the previous example shows.)
If this issue is resolved, the play leaves  $G^j$,
enters another connected component,
and the parties (players) proceed with the negotiation process.
Otherwise, the play cycles in  $G^j$  and negotiations terminate
resulting in the outcome  $c_j$.
We assume that it only matters at what issue
(a strongly connected component  $G^j$)
the parties have been forced into an impasse, while details
(a particular lasso, or a dicycle in which it ends in $G^j$)  are irrelevant.
Thus, an DGMS game is interpreted as a multistage diplomacy.
Then, the main result shows that the corresponding DGMS  game form is NS;
furthermore, for any  $u$  a NE in the game  $(g, u)$ can be found in linear time,
in other words, any bilateral multistage diplomatic conflict can be efficiently resolved
and the resolving strategies can be efficiently determined.

\section{Proofs}
\label{s1}
\subsection{Proof of Proposition \ref{p1}}
\label{ss11}


Let standardly  $I = \{1, \ldots, n\}$  and  $A = \{a_1, \ldots, a_p\}$
denote the sets of the players and outcomes, respectively.
We assume that  $n = 2$  and
prove that any two-person  win/lose-solvable game form
$g : X_1 \times X_2 \rightarrow A$  is NS.
The proof is constructive, for any utility function
$u : I \times A \rightarrow \R$  we will find a NE
in the obtained game  $(g, u)$  in at most  $2 p$  steps
in each of which we solve a  win/lose game
specified  by a  partition
$A = A_1 \cup A_2$,  where  $A_1$  and $A_2$  denote
the sets of outcomes winning for the players $1$ and $2$, respectively.

It will be convenient to replace  $u$  by two utility functions,
$u = (u_1, u_2), \;\; u_1 : A \rightarrow \R, \;  u_2 : A \rightarrow \R$,
of the first and second player, respectively.
Let us consider a partition  $A = W \cup W_1 \cup W_2$
satisfying the following properties:

\begin{itemize}
\item{(a1)}
any outcome of  $W_1$  is not better for player $1$
than any outcome of  $W$, that is,
$u_1(a') \leq u_1(a)$  for any  $a \in W$ and $a' \in W_1$;
\item{(a2)}
any outcome of  $W_2$  is not better for player $2$
than any outcome of  $W$, that is,
$u_2(a') \leq u_2(a)$  for any  $a \in W$ and $a' \in W_2$;
\item{(b1)}
player  $1$  cannot win (the win/lose game) with
$A_1 = W_2$, that is, $1$ cannot punish $2$;
\item{(b2)}
player  $2$  cannot win (the win/lose game) with
$A_2 = W_1$, that is, $2$  cannot punish $1$.
\end{itemize}

In general, it is not possible  that both players win, but
it is possible that both cannot win.
However, we assume that the game form  $g$
is win/lose-solvable, that is,
a player wins whenever the opponent cannot.


Our proof will be ``dynamical".
We begin with  $W = A$  and reduce  $W$
sending outcomes to  $W_1$  and to $W_2$
and verifying that all four above conditions always hold.
In the beginning, $W = A$, $W_1 = W_2 = \emptyset$, this is the case.
Let us show that  $W$  cannot be eliminated completely.
Indeed, suppose that
$W = \emptyset$  and, hence,  $A = W_1 \cup W_2$  is a partition.
Then, (b1), (b2), and  win/lose-solvability imply that
both players win the game with  $A_1 = W_2$  and  $A_2 = W_1$,
which is impossible.
Thus, reducing  $W$  must get stuck at some point, say, at a partition
$A = W \cup W_1 \cup W_2$.

Let  $a^* \in A$  be the worst outcome in $W$  for player $1$, that is,
$u_1(a^*) \leq u_1(a)$  for each $a \in W$.
We would like to move  $a^*$  from  $W$  to  $W_1$  but cannot.
The only possible reason is that (b2) will be violated,
that is, $2$  can win the game with $A_2 = W_1 \cup \{a^*\}$.
Denote by  $x_2$  a winning strategy.

Moving  $a^*$  alone from to  $W$ to $W_2$  may violate (a2).
To keep it, we have to add to  $a^*$  all
outcomes  $W_2(a^*) \subseteq W$  that are not better than
$a^*$  for $2$,  that is,
$u_2(a^*) \geq u_2(a)$  for all  $a \in W_2(a^*)$;
in particular, $a^* \in  W_2(a^*)$.
Then, moving $W_2(a^*)$  from  $W$  to $W_2$  is OK with (a2),
but still cannot be done.
The only possible reason is that (b1) will be violated,
that is,  $1$  can win the game with  $A_1 = W_2 \cup W_2(a^*))$.
Denote by  $x_1$  a winning strategy.

By construction,
$g(x_1, x_2) =  \{a^*\} = (W_1 \cup \{a^*\}) \cap (W_2 \cup W_2(a^*))$.
and  $(x_1, x_2)$  is a NE.

By each step, we are trying to send at least one outcome
either to  $W_1$  or to  $W_2$  and solve
the obtained  win/lose game to verify whether such sending can be realized.
Thus, totally, we will have to solve at most $2 p$  win/lose games.
\qed

\begin{remark}
The obtained NE is {\em simple}, that is,
$g(x_1) \cap g(x_2) = \{a^*\}$, where
$g(x_1) = \cup_{x_2} g(x_1, x_2)$  and
$g(x_2) = \cup_{x_1} g(x_1, x_2)$.
This property is important in applications;
see, for example, \cite{GK17}.
\end{remark}

\subsection{Proof of Proposition \ref{p2}}
\label{ss12}

Let us consider a two-person positional structure defined
by a digraph  $G = (V, E)$  and a partition  $V = V_1 \cup V_2 \cup V_T$.
Then, a  win/lose game is defined by an arbitrary partition
$A = \{c_j \mid j \in  J \setminus  J_0\} =  A_1 \cup A_2$
of the set of outcomes $A$  into the winning sets of players  $1$  and  $2$.

We suggest an algorithm that finds,
for a win/lose game, in time linear in $|V|$, a pair of
(pure stationary) strategies  $(x_1, x_2)$  that
form a subgame perfect saddle point.

\begin{remark}
Recall that no subgame  perfect NE may not exist when
the payoff function is not zero-sum; see the example of \cite[Section 6.2]{AGH10}.
Nevertheless, win/lose-solvability implies NS, although
the corresponding NE can be chosen subgame perfect only in the first case.
\end{remark}

We will combine solving zero-sum DG gasmes with backward induction.
Recall that  $G^*$  is a finite acyclic digraph
whose vertices $j \in J$  correspond
to the strongly connected components $G^j$  of  $G$.
Hence, $G^*$  must have a terminal vertex.
The corresponding  $G^j = (V^j, E^j)$  may be a terminal loop,
$j \in J_T$, but it may also be a non-trivial component,
with several vertices in it.
In this case, the same player wins
with respect to every initial position  $v_0 \in V^j$:
player $1$  wins if  $c_j \in A_1$  and  $2$ if  $c_j \in A_2$.
In both cases we eliminate all edges of  $E^j$
and add a terminal loop to each vertex of  $V^j$
that did not have any yet.
Thus, we can get rid of all terminal components $G^j$
such that  $j \in J \setminus (J_0 \cup J_T)$,
just increasing the number of terminal loops.

Although our algorithm modifies the original digraph  $G$
many times, yet, for simplicity, we will keep the notation.

Now let us choose a $j \in J \setminus (J_0 \cup J_T)$
such that component  $G^j = (V^j, E^j)$  is not terminal,
but any edge that leaves  $V^j$  enters a terminal loop.
More  precisely, for any edge  $e = (v, v') \in E$  such that  $v \in V^j$
there is a terminal loop  $e_{v'} = (v', v')$.
Note that loop  $e_v = (v, v)$  may exist but cannot be terminal.

We have to consider two cases: $c_j \in A_1$  and  $c_j \in A_2$.
Let us assume that
$c_j \in A_1$  and find the subset  $W \subseteq V^j$  such that
player  $2$  wins  whenever  $v_0 \in W$,
that is, (s)he can enforce a terminal loop from  $A_2$.
Set $W$ can be constructed by a standard recursive procedure.
We start with  $k = 0$  and set  $W^0$  of all terminal loops from  $A_2$.
Then,  for  $k = 0,1, \ldots $ we extend  $W^k$
adding to it each vertex  $v \in V^j \setminus W^k$  such that
\begin{itemize}
\item{(i)}
every move  $(v, v') \in E$  enters  $W^k$  if  $v \in V_1$, or
\item{(ii)}
there is a move  $(v, v') \in E$  that enters  $W^k$  if  $v \in V_2$.
\end{itemize}
We proceed until no vertex is added to  $W^{k^*}$  and then set  $W^{k^*} = W$.
By construction, player  $2$  wins whenever   $v_0 \in W$.
The winning strategy  $x_2$  is defined by the following rule:
stay in $W$  moving from  $W^k$ to $W^{k-1}$,
which is always possible, by (ii).
Then, sooner or later, the play will end in $W^0$.

Conversely, player $1$ wins whenever  $v_0 \in V^j \setminus W$.
The winning strategy  $x_1$  is defined by the following rule:
choose any move that does not enter $W$,  which is always possible, by (i).
Then, sooner or later, the play will realize either a terminal loop from $A_1$,
or a dicycle of  $G^j$. By our assumption, $c_j \in A_1$  too, hence,
player  $1$  wins in both cases.

If we assume that  $c_j \in A_2$ then the roles of players $1$ and $2$ swap,
otherwise the procedure remains the same.

\begin{remark}
The roles of the players differ in each case:
one who does not like  $c_j$  wins only if (s)he enforces
a winning terminal loop, while for the other one, who likes  $c_j$,
this is not necessary, it is enough just to stay in  $G^j$.
\end{remark}

The considered  win/lose game is solved whenever  $v_0 \in V^j$.
Eliminate all edges going from each vertex  $v \in V^j$ and
add a terminal loop  $e_v = (v, v)$  to it.
Set  $e_v \in A_2$  whenever  $v \in W$  and $e_v \in A_1$ otherwise.

In all above cases, we modified our win/lose game by eliminating  $G^j$  and
replacing it by terminal loops.
The number of the strongly connected components of  $G$
distinct from the terminal loops is reduced by one.
Then, we repeat the whole procedure.
In $|J| - |J_0|)$  steps the game will be solved
with respect to all  $v_0 \in V$.
\qed

\begin{remark}
By construction, the obtained NE in the considered
 win/lose game is subgame perfect.
In contrast, the original game  $(g,u)$  may
have no subgame perfect NE at all; see examples in
\cite[Figure 6]{AGH10} and \cite{BEGM12}.
\end{remark}

\subsection{Proof of Proposition \ref{p3}}
\label{ss13}


Consider a a two-person zero-sum DGMS game modeled by a digraph  $G$.
The proof follows the same line as the proof of Proposition \ref{p2}
until we start to analyze a component  $G^j = (V^j, E^j)$
from which every move comes to a terminal loop.
Obviously, this is a BD subgame and, according to \cite{AHMS12},
it has a subgame perfect saddle point than can be
determined in ``almost linear"  time  $(|E_j| + |V_j| + |A_j| \log |A_j|)$,
where  $A_j$  is the set of terminal loops reachable from  $V_j$.
These strategies defined on $V^j$  will be extended
by backward induction to the whole vertex-set  $V$.

As before, we eliminate all edges
beginning in  $V^j$  and add to each  $v \in V^j$
the terminal loop  $e_v$  assigning to it
the value obtained for the subgame.
Then, we choose the next component  $G^j$ and
repeat the whole procedure.
In $|J| - |J_0|)$  steps the game will be solved
with respect to all  $v_0 \in V$.
It is not difficult to see that, by construction,
we obtain a subgame perfect saddle point.
As for complexity, the ``almost linear" bound
of  \cite{AHMS12}  still holds, since for
$G^j$  the algorithm is the same and
the backward induction works  for all  $j \in J$  in linear time.
Thus, the worst case is  $|J \setminus (J_0 \cup J_T)| = 1$  corresponding to DG games.

\section{Open questions}
\label{s2}


As before, let  $n$  and  $p$  denote  the numbers of players and
terminals (or, what is the same, terminal loops), and let
$\ell$  be the length of the longest directed path in  $G^*$.
In the present paper we proved that NS holds for $n = 2$.

Yet, NS may fail to exist when $n > 2$.
A three-person NE free example is given in \cite[Figure 1]{BGMOV16}.
This example is pretty compact: $p = \ell = 3$;
in addition to three terminal loops, it contains only one dicycle
(of length $2$) and  two positions approaching it.
Note that each of them is a strongly connected component.
Only one player controls two positions, while other two control
a unique position each.

One can easily reduce $\ell =  3$  to  $2$
without creating a NE.
Consider the move $(1,2)$ in \cite[Figure 1]{BGMOV16} and
add to it the return  $(2,1)$.
Clearly, for the obtained digraph  $G'$ we have $\ell' = 2$.
Furthermore, by adding $(2,1)$  we create
the second dicycle  $c' = ((1,2),(2,1))$,
which becomes the fifth outcome, in addition to
terminals $a_1, a_2, a_3$, and dicycle  $c = ((3,2),(2,3))$.
We leave for the reader to verify that
no NE will appear, even if we identify  $c'$  with  $c$.

The example given in  \cite{Gur15}  was larger:
$n = \ell = 4$, $p = 5$, it also has a unique dicycle
(of length $5$) and three positions approaching it.

\medskip

It remains open whether NS results from any of the following four assumptions:
\begin{itemize}
\item{$(A1)$} $\ell = 1$;
\item{$(A2)$} $p \leq 2$;
\item{$(A3)$} each player controls a unique position
(so-called {\em play-once} case);
\item{$(A4)$} each terminal loop is better than any other dicycle for every player?
\end{itemize}

The last three questions were considered before, but only for the DG games.
Let us introduce such an assumption:

\medskip

$(A0)$  All dicycles form a single outcome  $c$.

\medskip

NS results from
(($(A2)$, $(A4)$, and $(A0)$)  or ($(A3)$, $(A4)$ and $(A0)$).
Both implications follow from the results of \cite{BG03}.
Note that assuming  $(A2)$ and  $(A4)$
we can also assume, without loss of generality, that  $n \leq 2$  and
then, NS follows, by Theorem \ref{t1}.

In the preprint \cite{BR09}, the authors replaced  $(A2)$, $p \leq 2$,
in ($(A2)$, $(A4)$, and $(A0)$) by a weaker assumption, $p \leq 3$, and derived NS.
However, they found an irreparable flaw in their proof
and, thus, the implication remains open.

\medskip

In theory of stochastic games
different types of {\em additive costs} are considered.
Let  $r : I \times E \rightarrow \R$  be a local cost function:
each player  $i \in I$  has to pay the amount  $r(i, e)$
whenever the directed edge  $e \in E$  is chosen as a move by a player.
Let us also assume

\smallskip

(A5) the local cost function $r$  takes only strictly positive values.

\smallskip

Note that $(A4)$ follows from  $(A5)$.
Now, we do not introduce terminal loops.
If a play  $P$  comes to a terminal,
we define the effective cost of  $P$  for each player  $i \in I$
as the sum of all local costs for all edges of the play,
that is, $U(i, P) = \sum_{e \in P} r(i, e)$.
If $P$  is a lasso then  $U(i, P) = + \infty$  for all  $i \in I$.
All players minimize their effective costs.

The main result of \cite{BG03} states that
in this case NS follows from $(A3)$.

\begin{conjecture}
In the two-person case, $(A5)$  implies  NS.
\end{conjecture}

A NE free example for  $n = 3$,
and for  $n=2$  with $(A5)$  waved, were constructed in \cite{GO14}.

\medskip

{\bf Acknowledgements:}
This work was partially funded by Russian Academic Excellence Project '5-100'.
The author is thankful to Endre Boros for helpful remarks and suggestions.


\begin{thebibliography}{99}

\bibitem{AGH10}
D. Anderson, V. Gurvich, and T. Hansen,
On acyclicity of games with cycles,
Discrete Applied Math. 158 (10) (2010) 1049-1063.

\bibitem{AHMS12}
D. Andersson, K. Hansen, P. Miltersen, and T. Sorensen,
Deterministic graphical games, revisited,
J. Logic and Computation 22:2 (2012) 165-178.
Preliminary version in Fourth Conference on Computability in Europe (CiE-08),
Lecture Notes in Computer Science 5028 (2008) 1--10.

\bibitem{BBM89}
A.B. Barabas, R.E. Basko, and I.S. Menshikov,
On an approach to analysis of conflict situations,
Proceedings of the  Moscow Computing Center 2 (1989) 13--20.

\bibitem{BG03}
E. Boros and V. Gurvich,
On Nash-solvability in pure strategies of finite games
with perfect information which may have cycles,
Math. Soc. Sciences 46 (2003) 207-241.

\bibitem{BG09}
E. Boros and V. Gurvich,
Why Chess and Backgammon can be solved in pure
positional uniformly optimal strategies,
RUTCOR Research Report, RRR-21-2009, Rutgers University.

\bibitem{BGMS07}
E. Boros, V. Gurvich, K. Makino, and Wei Shao,
Nash-solvable two-person symmetric cycle game forms,
Discrete Applied Math. 159 (2011) 1461--1487.

\bibitem{BEGM12}
E. Boros, K. Elbassioni, V. Gurvich, and K. Makino.
On Nash Equilibria and Improvement Cycles in Pure Positional Strategies
for Chess-like and Backgammon-like n-person Games,
Discrete Math. 312:4 (2012) 772--788.

\bibitem{BGMOV16}
E. Boros, V. Gurvich, M. Milanic, V. Oudalov, and J. Vicic,
A three-person deterministic graphical game without Nash equilibria,
submitted to Discrete Applied Mathematics, October 16, 20016;
now available at: https://scirate.com/arxiv/1610.07701


\bibitem{BR09}
E. Boros and R. Rand,
Terminal games with three terminals have proper Nash equilibria,
RUTCOR Research Report, RRR-22-2009, Rutgers University.

\bibitem{Con92}
A. Condon,
The complexity of stochastic games,
Information and Computation 96 (1992) 203--224.


\bibitem{DS91}
V.I. Danilov, A.I. Sotskov, Social Choice Mechanisms,
Moscow Nauka 1991 (in Russian);
English translation in ``Studies of Economic Design",
Springer, Berlin-Heidelberg, 2002.

\bibitem{EF70}
J. Edmonds and D.R. Fulkerson, Bottleneck extrema,
J. of Combinatorial Theory, 8 (1970),  299-306.

\bibitem{Gal53}
D. Gale, A theory of $N$-person games with perfect information,
Proc. Natl. Acad. Sci. 39 (1953) 496--501.

\bibitem{Gil57}
D. Gillette,
Stochastic games with zero stop probabilities,
Contribution to the Theory of Games III",
Annals of Mathematics Studies 39 (1957) 179--187.


%
%

\bibitem{Gur73}
V. Gurvich, To theory of multi-step games,
USSR Comput. Math. and Math. Phys. 13:6 (1973) 143-161.

\bibitem{Gur75}
V. Gurvich, Solution of positional games in pure strategies,
USSR Comput. Math. and Math. Phys. 15:2 (1975) 74-87.


\bibitem{Gur88}
V. Gurvich, Equilibrium in pure strategies,
Soviet Math. Dokl. 38:3 (1989) 597-602.

\bibitem{Gur97}
V. Gurvich, Criteria of nonemptiness for dual cores,
Doklady Mathematics  55:1 (1997) 12--16.


\bibitem{Gur15}
V. Gurvich. 
A four-person chess-like game
without Nash equilibria in pure stationary strategies,
Business Informatics 1:31 (2015) 68--76,
arXiv 1411.0349   http://arxiv.org/abs/1411.0349 .

\bibitem{GK17}
V. Gurvich and G. Koshevoy.
Monotone bargaining is Nash-solvable;
submitted to Mathematics of Operations Research 11/4/2017, ID MOR-2017-309;
available online at  ArXiv 11/3/2017  http://arxiv.org/abs/1711.00940.

\bibitem{GO14}
V. Gurvich and V. Oudalov.
On Nash-solvability in pure stationary strategies
of the deterministic n-person games
with perfect information and mean or total effective cost,
Discrete Appl. Math. 167 (2014) 131--143.

\bibitem{HL31}
G. H. Hardy and J. E. Littlewood,
Notes on the theory of series (XVI): two Tauberian theorems,
J. of London Mathematical Society 6 (1931) 281--286.

\bibitem{Kuh50}
H. Kuhn, Extensive games,
Proc. Natl. Acad. Sci. 36 (1950) 286--295.

\bibitem{Kuh53}
H. Kuhn, Extensive games and the problem of information,
in Contributions to the theory of games, Volume 2, Princeton (1953) 193--216.

\bibitem{LL69}
T. M. Liggett and S. A. Lippman,
Stochastic games with perfect information and time-average payoff,
SIAM Review 4 (1969) 604--607.


\bibitem{Nas50}
J. Nash, Equilibrium points in n-person games,
Proceedings of the National Academy of Sciences 36:1 (1950) 48--49.

\bibitem{Nas51}
J. Nash, Non-Cooperative Games,
Annals of Mathematics 54:2 (1951) 286--295.




\bibitem{RM-DS17}
S. Le Roux, E. Martin-Dorel, and J.G. Smaus,
An existence theorem of Nash equilibrium in Coq and Isabelle,
Proceedings of the Eighth International Symposium on Games,
Automata, Logics, and Formal Verification
(GandALF) Rome, Italy,  Septembber 20--22 (2017) 46--60.

\bibitem{Sha81}
M. Sharir, A strong-connectivity algorithm and its application in data flow analysis,
Comput. Math. Appl. 7 (1981) 67--72.

\bibitem{Tar72}
R.E. Tarjan, Depth-first search and linear graph algorithms,
SIAM J. Computing 1:2 (1972) 146--160.
doi:10.1137/0201010

\bibitem{Tho86}
K. Thompson, Retrograde analysis of certain endgames,
Journal of the International Computer Chess Association
9:3 (1986) 131--139.

\bibitem{Was90}
A. R. Washburn, Deterministic graphical games,
J. of Math. Analysis and Applications 153 (1990) 84--96.

\end{thebibliography}
\end{document}